\documentclass[11pt, leqno]{article}
\usepackage{euscript,amstext,amsmath,amsthm,a4,amssymb}

\newcommand{\tr}{{\rm Tr}}

\newcommand{\comment}[1]{}

\newtheorem{thm}{Theorem}[section]
\newtheorem{prop}[thm]{Proposition}
 \newtheorem{cor}[thm]{Corollary}

\theoremstyle{remark}
\newtheorem{Rem}[thm]{Remark}
\theoremstyle{definition}

\title{$\eta$-invariant  and flat vector bundles}
\author{Xiaonan Ma\footnote{Centre de Math\'ematiques, UMR 7640 du CNRS,
\'Ecole Polytechnique, 91128 Palaiseau Cedex, France.
(ma@math.polytechnique.fr)}\ \ and\ \ Weiping
Zhang\footnote{Nankai Institute of Mathematics \& LPMC, Nankai
University, Tianjin 300071, P.R. China. (weiping@nankai.edu.cn)}}
\date{\it Dedicated to the memory of Professor Shiing-shen Chern}

\begin{document}

\maketitle
\begin{abstract} We present an alternate definition of the mod
{\bf Z} component of the Atiyah-Patodi-Singer $\eta$ invariant
associated to (not necessary unitary) flat vector bundles, which
identifies explicitly its  real and imaginary parts.
 This is
done by combining a deformation of flat connections introduced in
a previous paper with the analytic continuation procedure
appearing  in the original article of Atiyah, Patodi and Singer.

$\ $


{\bf 2000 MR Subject Classification} 58J
\end{abstract}

\renewcommand{\theequation}{\thesection.\arabic{equation}}
\setcounter{equation}{0}

\section{Introduction} \label{s1}

Let $M$ be an odd dimensional oriented closed  spin manifold
carrying a Riemannian metric $g^{TM}$. Let $S(TM)$ be the
associated Hermitian bundle of spinors. Let $E$ be a Hermitian
vector bundle over $M$ carrying a unitary connection $\nabla^E$.
Moreover, let $F$ be a Hermitian vector bundle over $M$ carrying a
unitary flat connection $\nabla^F$. Let
\begin{align}\label{1.1} D^{E\otimes F}:\Gamma(S(TM)\otimes E\otimes F)\longrightarrow
\Gamma(S(TM)\otimes E\otimes F)
\end{align}
denote the corresponding (twisted) Dirac operator, which is
formally self-adjoint (cf. \cite{BGV}).

For any $s\in {\bf C}$ with ${\rm Re}(s)>>0$, following
\cite{APS1}, set
\begin{align}\label{1.2}
\eta\left(D^{E\otimes F},s\right)=\sum_{\lambda\in{\rm
Spec}(D^{E\otimes F})\setminus \{ 0\}}{{\rm Sgn}(\lambda)\over
|\lambda|^s}.
\end{align} Then by \cite{APS1}, one knows that $\eta (D^{E\otimes F},s
)$ is a holomorphic function in $s$ when ${\rm Re}(s)>{\dim M\over
2}$. Moreover, it extends to  a meromorphic function over ${\bf
C}$, which is holomorphic at $s=0$. The $\eta$ invariant of
$D^{E\otimes F}$, in the sense of Atiyah-Patodi-Singer
\cite{APS1}, is defined by
\begin{align}\label{1.3}\eta\left(D^{E\otimes F} \right)=
\eta\left(D^{E\otimes F},0\right) ,
\end{align}
while the corresponding {\it reduced} $\eta$ invariant is defined
and denoted  by
\begin{align}\label{1.4}\overline{\eta}\left(D^{E\otimes F} \right)=
{\dim \left(\ker D^{E\otimes F}\right)+\eta\left(D^{E\otimes
F}\right)\over 2} .
\end{align}

The $\eta$ and reduced $\eta$ invariants play an important role in
the Atiyah-Patodi-Singer index theorem for Dirac operators on
manifolds with boundary \cite{APS1}.

In \cite{APS2} and \cite{APS3}, it is shown that the following
quantity
\begin{align}\label{1.5}\rho\left(D^{E\otimes F}\right):=\overline{\eta}\left(D^{E\otimes F}
\right)-{\rm rk}(F)\, \overline{\eta}\left(D^{ E} \right)\ \ \
{\rm mod}\ \ {\bf Z}
\end{align}
does not depend on the choice of $g^{TM}$ as well as the metrics
and (Hermitian) connections on $E$. Also, a Riemann-Roch theorem
  is proved in
\cite[(5.3)]{APS3}, which gives a $K$-theoretic  interpretation of
the analytically defined invariant $\rho (D^{E\otimes F} )\in{\bf
R/Z }$. Moreover, it is pointed out in \cite[Page 89, Remark
(1)]{APS3} that the above mentioned $K$-theoretic interpretation
applies also to the case where $F$ is a non-unitary flat vector
bundle, while on \cite[Page 93]{APS3} it shows how one can define
the reduced $\eta$-invariant in case $F$ is non-unitary, by
working on non-self-adjoint elliptic operators, and then extend
the Riemann-Roch result \cite[(5.3)]{APS3} to an identity in {\bf
C/Z} (instead of {\bf R/Z}). The idea of analytic continuation
plays a key role in obtaining this Riemann-Roch result, as well
as its non-unitary extension.

In this paper, we show that by using the idea of analytic
continuation, one can construct the {\bf C/Z} component of
$\overline{\eta} (D^{E\otimes F} )$ directly, with out passing to
analysis of non-self-adjoint operators, in case where $F$ is a
non-unitary flat vector bundle. Consequently, this leads to a
direct construction of $\rho (D^{E\otimes F} )$ in this case. We
will   use a deformation introduced in \cite{MZ} for flat
connections in our construction.

In the next section, we will first recall the above mentioned
deformation from \cite{MZ} and then give our construction of
$\overline{\eta} (D^{E\otimes F} )$ mod {\bf Z} and $\rho
(D^{E\otimes F} )\in {\bf C/Z}$ in the case where $F$ is a
non-unitary flat vector bundle.

$\  $

\noindent {\bf Acknowledgements}. The work of the second author
was partially supported by the Cheung-Kong Scholarship of the
Ministry of Education of China and the 973 Project of the Ministry
of Science and Technology of China.

\section{The $\eta$ and $\rho$ invariants associated to non-unitary flat
vector bundles}\label{s2}

This section is organized as follows. In Section \ref{2a}, we
construct certain secondary characteristic forms and classes
associated to non-unitary flat vector bundles. In Section
\ref{2b}, we present our construction of the mod {\bf Z} component
of the reduced $\eta$-invariant, as well as the $\rho$-invariant,
associated to non-unitary flat vector bundles. Finally, we include
some further remarks in Section \ref{2c}.

\subsection{Chern-Simons classes and flat vector
bundles}\label{2a} \setcounter{equation}{0} We fix a square root
of $\sqrt{-1}$ and let
$\varphi:\Lambda(T^*M)\rightarrow\Lambda(T^*M)$ be the
homomorphism defined by
$\varphi:\omega\in\Lambda^i(T^*M)\rightarrow (2\pi
\sqrt{-1})^{-i/2}\omega.$ The formulas in what follows will not
depend on the choice of the square root of $\sqrt{-1}$.

If  $W$ is a  complex vector bundles over $M$ and  $\nabla^W_0$,
$\nabla^W_1$ are two connections on $W$. Let $W_t$, $0\leq t\leq
1$, be a smooth path of connections on $W$ connecting $\nabla^W_0$
and $\nabla^W_1$. We define Chern-Simons form
$CS(\nabla^W_0,\nabla^W_1)$ to be the differential form given by
\begin{align}\label{2.5}
CS(\nabla^W_0,\nabla^W_1)=-\left(1\over
2\pi\sqrt{-1}\right)^{1\over 2}
\varphi\int_0^1\tr\left[{\partial\nabla^W_t\over\partial
t}\exp(-(\nabla^W_t)^2) \right]dt.
\end{align}
Then (cf. \cite[Chapter 1]{Z1})
\begin{align}\label{2.6}
dCS(\nabla^W_0,\nabla^W_1)={\rm ch}(W,\nabla^W_1)-{\rm
ch}(W,\nabla^W_0).
\end{align}
Moreover, it is well-known that up to exact forms,
$CS(\nabla^W_0,\nabla^W_1)$ does not depend on the path of
connections on $W$ connecting $\nabla^W_0$ and $\nabla^W_1$.

 Let $(F,\nabla^F)$ be
a flat vector bundle carrying the flat connection $\nabla^F$. Let
$g^F$ be a Hermitian metric on $F$. We do not assume that
$\nabla^F$ preserves $g^F$. Let $(\nabla^F)^*$ be the adjoint
connection of $\nabla^F$ with respect to $g^F$.

From \cite[(4.1), (4.2)]{BZ} and \cite[\S 1(g)]{BL}, one has
\begin{align}\label{2.1}
\left(\nabla^F\right)^*=\nabla^F+\omega\left(F,g^F\right)
\end{align}
with
\begin{align}\label{2.2}
 \omega\left(F,g^F\right)=\left(g^F\right)^{-1}\left(\nabla^Fg^F\right).
\end{align} Then
\begin{align}\label{2.3}
  \nabla^{F,e}=\nabla^F+{1\over 2}\omega\left(F,g^F\right)
\end{align}
is a Hermitian connection on $(F,g^F)$ (cf. \cite[(1.33)]{BL} and
\cite[(4.3)]{BZ}).

Following \cite[(2.47)]{MZ}, for any $r\in{\bf C}$, set
\begin{align}\label{2.4}
  \nabla^{F,e,(r)}=\nabla^{F,e}+{\sqrt{-1}r\over
  2}\omega\left(F,g^F\right).
\end{align}
Then for any $r\in{\bf R}$, $\nabla^{F,e,(r)}$ is a Hermitian
connection on $(F,g^F)$.

On the other hand, following \cite[(0.2)]{BL},  for any integer
$j\geq 0$, let $c_{2j+1}(F,g^F)$ be the Chern form defined by
\begin{align}\label{2.7}
c_{2j+1}\left(F,g^F\right)=(2\pi\sqrt{-1})^{-j}2^{-(2j+1)}\tr\left[\omega^{2j+1}\left(F,g^F\right)\right].
\end{align} Then $c_{2j+1}(F,g^F)$ is a closed form on $M$.
Let $c_{2j+1}(F)$ be the associated cohomology class in
$H^{2j+1}(M,{\bf R})$, which does not depend on the choice of
$g^F$.

For any $j\geq 0$ and $r\in {\bf R}$, let $a_j(r)\in {\bf R}$ be
defined as
\begin{align}\label{2.8}
a_j(r)=\int_0^1\left(1+u^2r^2\right)^jdu.
\end{align}

With these notation we can now state the following result first
proved in \cite[Lemma 2.12]{MZ}.

\begin{prop}\label{t2.1} The following identity in
$H^{\rm odd}(M,{\bf R})$ holds for any  $r\in {\bf R}$,
\begin{align}\label{2.9}
CS\left(\nabla^{F,e},\nabla^{F,e,(r)}\right)=-{r\over
2\pi}\sum_{j=0}^{+\infty} {a_j(r)\over j!}c_{2j+1}(F).
\end{align}
\end{prop}

\subsection{$\eta$ and $\rho$ invariants associated to flat vector
bundles}\label{2b}

We now make the same assumptions as in the beginning of Section 1,
except  that we no longer assume $\nabla^F$ there is unitary.

For any $r\in {\bf C}$,  let
\begin{align}\label{2.10}D^{E\otimes F}(r):\Gamma(S(TM)\otimes E\otimes F)\longrightarrow
\Gamma(S(TM)\otimes E\otimes F)
\end{align}
 denote the Dirac operator associated to the  connection
 $\nabla^{F,e,(r)}$ on $F$. Since when $r\in{\bf R}$, $\nabla^{F,e,(r)}$ is Hermitian on $(F,g^F)$, $D^{E\otimes
 F}(r)$ is formally self-adjoint and one can define the associated
 reduced $\eta$-invariant as in (\ref{1.4}).

 By the variation formula for the
 reduced
 $\eta$-invariant (cf. \cite{APS1} and \cite{BF}), one gets that for any $r\in {\bf R}$,
\begin{align}\label{2.11}\overline{\eta}\left(D^{E\otimes
F}(r)\right)-\overline{\eta}\left(D^{E\otimes F}(0)\right)\equiv
\int_M\widehat{A}(TM){\rm
ch}(E)CS\left(\nabla^{F,e},\nabla^{F,e,(r)}\right)\ \ {\rm mod}\
{\bf Z},
\end{align}
where $\widehat{A}$ and ${\rm ch}$ are standard notations for the
Hirzebruch $\widehat{A}$-class and Chern character respectively
(cf. \cite[Chapter 1]{Z1}).

Let $D^{E\otimes F,e}$ denote the Dirac operator $D^{E\otimes
F}(0)$.

From (\ref{2.9}) and (\ref{2.11}), one gets that for any $r\in{\bf
R}$,
\begin{align}\label{2.12}\overline{\eta}\left(D^{E\otimes
F}(r)\right)\equiv \overline{\eta}\left(D^{E\otimes
F,e}\right)-{r\over 2\pi}\int_M\widehat{A}(TM){\rm ch}(E)
\sum_{j=0}^{+\infty} {a_j(r)\over j!}c_{2j+1}(F)\ \ {\rm mod}\
{\bf Z}.
\end{align}

Recall that even though when ${\rm Im}(r)\neq 0$, $D^{E\otimes
F}(r)$ might not be formally self-adjoint, the $\eta$-invariant
can still be defined, as outlined in \cite[page 93]{APS3}. On the
other hand, from (\ref{2.3}) and (\ref{2.4}), one sees that
\begin{align}\label{2.13}\nabla^F=\nabla^{F,e,(\sqrt{-1})}.
\end{align}
We denote the associated Dirac operator $D^{E\otimes
F}(\sqrt{-1})$ by $D^{E\otimes F}$.

We also recall that
\begin{align}\label{2.14.1}
 \int_0^1\left(1-u^2\right)^jdu={2^{2j}(j!)^2\over (2j +1)!}.
\end{align}

We can now state the main result of this paper as follows.

\begin{thm} \label{t2.2} Formula (\ref{2.12}) holds indeed for any
$r\in{\bf C}$. In particular, one has
\begin{align}\label{2.14}
\overline{\eta}\left(D^{E\otimes F}\right)\equiv
\overline{\eta}\left(D^{E\otimes F,e}\right)-{\sqrt{-1}\over
2\pi}\int_M\widehat{A}(TM){\rm ch}(E) \sum_{j=0}^{+\infty} {
2^{2j}j!\over (2j+1)!}c_{2j+1}(F)\ \ {\rm mod}\ {\bf Z}.
\end{align}
Equivalently,
\begin{align}\label{a2.14}
&{\rm Re} \left(\overline{\eta}\left(D^{E\otimes F}\right)\right)
\equiv \overline{\eta}\left(D^{E\otimes F,e}\right) \ \ {\rm mod}\ {\bf Z},\\
&{\rm Im} \left(\overline{\eta}\left(D^{E\otimes F}\right)\right)
=-\frac{1}{2\pi}\int_M\widehat{A}(TM){\rm ch}(E)
\sum_{j=0}^{+\infty} { 2^{2j}j!\over (2j+1)!}c_{2j+1}(F).
\nonumber
\end{align}
\end{thm}
{\it Proof}. Clearly, the right hand side of (\ref{2.12}) is a
holomorphic function in $r\in{\bf C}$. On the other hand, by
\cite[page 93]{APS3}, $\overline{\eta} (D^{E\otimes F}(r) )$ mod
{\bf Z} is also holomorphic in $r\in{\bf C}$. By (\ref{2.12}) and
the uniqueness of the analytic continuation, one sees that
(\ref{2.12}) holds indeed for any $r\in{\bf C}$. In particular,
by putting together (\ref{2.12}) and (\ref{2.13}), one gets
(\ref{2.14}).\ \ \ \ Q.E.D.

$\ $

Recall that when $\nabla^F$ preserves $g^F$, the $\rho$-invariant
has been defined in (\ref{1.5}). Now if we no longer assume that
$\nabla^F$ preserves $g^F$, then by Theorem \ref{t2.2}, one sees
that one gets the following formula of the associated  (extended)
$\rho$-invariant.

\begin{cor}\label{t2.3} The following identity holds,
\begin{multline}\label{2.15}
\rho\left(D^{E\otimes
F}\right)\equiv\overline{\eta}\left(D^{E\otimes F,e} \right)-{\rm
rk}(F)\, \overline{\eta}\left(D^{ E} \right)\\
-{\sqrt{-1}\over 2\pi}\int_M\widehat{A}(TM){\rm ch}(E)
\sum_{j=0}^{+\infty} { 2^{2j}j!\over (2j+1)!} c_{2j+1}(F)\ \ \
{\rm mod}\ \ {\bf Z}.
\end{multline}
Equivalently,
\begin{align}\label{a2.14c}
&{\rm Re} \left(\rho\left(D^{E\otimes F}\right)\right) \equiv
\overline{\eta}\left(D^{E\otimes F,e}\right)- {\rm
rk}(F)\, \overline{\eta}\left(D^{ E} \right)\ \ {\rm mod}\ {\bf Z},\\
&{\rm Im} \left(\rho\left(D^{E\otimes F}\right)\right)
=-\frac{1}{2\pi}\int_M\widehat{A}(TM){\rm ch}(E)
\sum_{j=0}^{+\infty} { 2^{2j}j!\over (2j+1)!}c_{2j+1}(F).
\nonumber
\end{align}
\end{cor}

It is pointed out in \cite{APS3} that the Riemann-Roch formula
proved in \cite[(5.3)]{APS3} still holds for $\rho (D^{E\otimes F}
)$ in the case where $\nabla^F$ does not preserve $g^F$. One way
to understand this is that the argument in the proof of
\cite[(5.3)]{APS3} given in \cite{APS3} works line by line to give
a $K$-theoretic interpretation of $\overline{\eta} (D^{E\otimes
F,e}  )-{\rm rk}(F)\, \overline{\eta} (D^{ E}  )$. By (\ref{2.15})
it then gives such an interpretation for $\rho (D^{E\otimes F} )$.

\subsection{Further remarks}\label{2c}

\begin{Rem}\label{t2.21}  The argument in proving Theorem
\ref{t2.2} works indeed for any twisted vector bundles $F$, not
necessary a flat vector bundle. This gives a direct formula for
the mod {\bf Z} part of the $\eta$-invariant for non-self-adjoint
Dirac operators.
\end{Rem}

\begin{Rem}
In \cite[Theorem 2.2]{Z2}, a $K$-theoretic formula for
$D^{E\otimes F}(r)$ mod {\bf Z} has been given in the $r\in{\bf
R}$ case.   As a consequence, one gets an alternate $K$-theoretic
formula for $\rho (D^{E\otimes F} )$  in \cite[(4.6)]{Z2} which
holds  in the case where $\nabla^F$   preserves  $g^F$. By
combining the arguments in \cite{Z2} with Theorem \ref{t2.2}
proved above, one can indeed extend \cite[Theorem 2.2]{Z2} and
\cite[(4.6)]{Z2} to the case where $\nabla^F$  might not preserve
$g^F$. We leave this to the interested reader. Here we only
mention that this will provide an alternate $K$-theoretic
interpretation of   $\rho$-invariants in the case where $\nabla^F$
does not preserve $g^F$.
\end{Rem}

\begin{Rem}
We refer to \cite{MZ} where we have employed  deformation
(\ref{2.4}) to study and generalize certain
Riemann-Roch-Grothendieck formulas due to Bismut-Lott (\cite{BL})
and Bismut (\cite{B1}), for flat vector bundles over fibred
spaces.
\end{Rem}

\end{document}